
\magnification1200
\input amstex.tex
\documentstyle{amsppt}

\hsize=12.5cm \vsize=18cm
\hoffset=1cm \voffset=2cm

\footline={\hss{\vbox to 2cm{\vfil\hbox{\rm\folio}}}\hss}
\nopagenumbers
\def\DJ{\leavevmode\setbox0=\hbox{D}\kern0pt\rlap
{\kern.04em\raise.188\ht0\hbox{-}}D}

\def\txt#1{{\textstyle{#1}}}
\baselineskip=13pt
\def\hf{{\textstyle{1\over2}}}
\def\a{\alpha}\def\b{\beta}
\def\d{{\,\roman d}}
\def\e{\varepsilon}
\def\f{\varphi}
\def\G{\Gamma}
\def\k{\kappa}
\def\s{\sigma}

\def\={\;=\;}

\def\zt{\zeta(\hf+it)}

\def\R{\Re{\roman e}\,} \def\I{\Im{\roman m}\,}
\def\z{\zeta}

\def\H{H_j^4({\txt{1\over2}})} 
\def\hf{{\textstyle{1\over2}}}
\def\txt#1{{\textstyle{#1}}}
\def\f{\varphi}

\font\tenmsb=msbm10
\font\sevenmsb=msbm7
\font\fivemsb=msbm5
\newfam\msbfam
\textfont\msbfam=\tenmsb
\scriptfont\msbfam=\sevenmsb
\scriptscriptfont\msbfam=\fivemsb
\def\Bbb#1{{\fam\msbfam #1}}

\def \NN {\Bbb N}

\def \RR {\Bbb R}
\def \ZZ {\Bbb Z}

\font\ff=cmr8
\def\txt#1{{\textstyle{#1}}}
\baselineskip=13pt

\font\teneufm=eufm10
\font\seveneufm=eufm7
\font\fiveeufm=eufm5
\newfam\eufmfam
\textfont\eufmfam=\teneufm
\scriptfont\eufmfam=\seveneufm
\scriptscriptfont\eufmfam=\fiveeufm
\def\mathfrak#1{{\fam\eufmfam\relax#1}}

\font\tenmsb=msbm10
\font\sevenmsb=msbm7
\font\fivemsb=msbm5
\newfam\msbfam
     \textfont\msbfam=\tenmsb
      \scriptfont\msbfam=\sevenmsb
      \scriptscriptfont\msbfam=\fivemsb
\def\Bbb#1{{\fam\msbfam #1}}

\def \NN {\Bbb N}

\def \RR {\Bbb R}
\def \ZZ {\Bbb Z}

  \def\rightheadline{{\hfil{\ff
  On the moments of Hecke series at central points}\hfil\tenrm\folio}}

  \def\leftheadline{{\tenrm\folio\hfil{\ff
   A. Ivi\'c }\hfil}}
  \def\emptyheadline{\hfil}
  \headline{\ifnum\pageno=1 \emptyheadline\else
  \ifodd\pageno \rightheadline \else \leftheadline\fi\fi}

\topmatter
\title
ON THE MOMENTS OF HECKE SERIES AT CENTRAL POINTS
\endtitle
\author   Aleksandar Ivi\'c  \endauthor
\address
Aleksandar Ivi\'c, Katedra Matematike RGF-a
Universiteta u Beogradu, \DJ u\v sina 7, 11000 Beograd,
Serbia
\endaddress
\email {\tt
ivic\@rgf.bg.ac.yu, aivic\@matf.bg.ac.yu} \endemail
\dedicatory
Functiones et Approximatio XXX(2002), 7-40
\enddedicatory
\abstract
{\ff We prove, in standard notation from spectral theory,
the following asymptotic formulas:
$$
\sum_{\k_j\le K} \a_jH_j^3(\hf) = K^2P_3(\log K) + O(K^{5/4}\log^{37/4}K)
$$
and
$$
\sum_{\k_j\le K} \a_j\H = K^2P_6(\log K) + O(K^{3/2}\log^{25/2}K),
$$
where $P_3(x)$ and $P_6(x)$ are polynomials of degree three and six,
whose coefficients may be explicitly evaluated.}
\endabstract
\endtopmatter

\document

\head
1. Introduction and statement of results
\endhead

The purpose of this paper is to obtain asymptotic formulas for
sums of $H_j^3(\hf)$ and $\H$, where $H_j(s)$ is the Hecke series,
to be defined below. Sums with $H_j(\hf)$ are important for several
reasons, one of
which is that they appear in the spectral decomposition of weighted
integrals involving $|\zt|^4$, which is of fundamental importance
in the theory of the Riemann zeta-function $\z(s)$.

\smallskip
We shall first present the
relevant notation  involving the spectral theory of
the non-Euclidean Laplacian will be given below.
For a competent and extensive account of
spectral theory the reader is referred to Y. Motohashi's monograph [15].

\medskip
Let $\,\{\lambda_j = \kappa_j^2 + {1\over4}\} \,\cup\, \{0\}\,$ be the
eigenvalues (discrete spectrum) of the hyperbolic Laplacian
$$
\Delta=-y^2\left({\left({\partial\over\partial x}\right)}^2 +
{\left({\partial\over\partial y}\right)}^2\right)
$$
acting over the Hilbert space composed of all
$\Gamma$-automorphic functions which are square integrable with
respect to the hyperbolic measure  ($\Gamma = \roman {PSL}(2,\ZZ)$).
Let $\{\psi_j\}_{j=1}^\infty$ be a maximal orthonormal system such that
$\Delta\psi_j=\lambda_j\psi_j$ for each $j\ge1$ and
$T(n)\psi_j=t_j(n)\psi_j$  for each integer $n\in\NN$, where
$$
\bigl(T(n)f\bigr)(z)\;=\;{1\over\sqrt{n}}\sum_{ad=n}
\,\sum_{b=1}^df\left({az+b\over d}\right)
$$
is the Hecke operator. We shall further assume that
$\psi_j(-\bar{z})=\e_j\psi_j(z)$ with $\e_j=\pm1$. We
then define ($s = \s + it$ will denote a complex variable)
$$
H_j(s)\;=\;\sum_{n=1}^\infty t_j(n)n^{-s}\qquad(\s > 1),
$$
which is the Hecke series associated with the Maass wave form
$\psi_j(z)$, and which can be continued to an entire function.
It  satisfies the functional equation
$$
H_j(s) = 2^{2s-1}\pi^{2s-2}\G(1-s+i\k_j)\G(1-s-i\k_j)
(\e_j\cosh(\pi\k_j)-\cos(\pi s))H_j(1-s),
$$
which by the Phragm\'en--Lindel\"of principle (convexity) implies the bound
$$
H_j(\hf) \;\ll_\e \; \k_j^{{1\over2}+\e},\leqno(1.1)
$$
where here and later $\e$ denotes arbitrarily small, positive  constants, not
necessarily the same ones at each occurrence. It is also important to note
that, from the work of Katok--Sarnak [9], it is known that $H_j(\hf) \ge 0$.

\smallskip
The sharpest asymptotic formula for sums of
$\a_jH^2_j(\hf)$ is due to Y. Motohashi [14]. His result is
$$
\sum_{\k_j\le T}\a_jH^2_j(\hf) = 2\pi^{-2}T^2(\log T + \gamma - \hf
- \log(2\pi)) + O(T\log^6T),\leqno(1.2)
$$
where $\gamma$ is Euler's constant,
$$
\a_j = |\rho_j(1)|^2(\cosh\pi\kappa_j)^{-1},
$$
and $\rho_j(1)$ is the first Fourier coefficient of  $\psi_j(z)$.

In what concerns known results
on sums of $\a_jH_j^3(\hf)$ and $\a_jH_j^4(\hf)$ we have
(see [15, Chapter 3])
$$
\sum_{\k_j\le K}\a_jH_j^4(\hf)  \ll K^2\log^{15}K\leqno(1.3)
$$
and
$$
\sum_{j=1}^\infty \a_jH_j^3(\hf)h_0(\k_j) =
\left({8\over3} + O\left({1\over\log K}\right)\right)\pi^{-3/2} K^3G
\log^3K\leqno(1.4)
$$
with
$$
K^{{1\over2}}\log^5K \le G \le K^{1-\e}, \leqno(1.5)
$$
$$
h_0(r)
= (r^2 + {\txt{1\over4}})\left\{\exp\left(-\left({r-K\over G}\right)^2\right)
+ \exp\left(-\left({r+K\over G}\right)^2\right)\right\}.\leqno(1.6)
$$
In [5] the author proved that
$$
\sum_{K\le\k_j\le K+1}\a_j H_j^3(\hf) \ll_\e K^{1+\e}.\leqno(1.7)
$$
This result immediately implies, since $H_j(\hf) \ge 0$ and
$\a_j \gg \k_j^{-\e}$ (see H. Iwaniec [6]), that
$$
H_j(\hf) \ll_\e \k_j^{{1\over3}+\e}, \leqno(1.8)
$$
which improves the convexity bound (1.1), and represents
hitherto the sharpest known unconditional upper bound for
$H_j(\hf)$. The bound (1.8) also follows from the result of
M. Jutila [7], namely
$$
\sum_{K\le\k_j\le K+K^{1/3}}\a_j\H \ll_\e K^{{4\over3}+\e},\leqno(1.9)
$$
and an extension of the bound (1.9) to sums of $|H_j(\hf + it)|^4$
has been attained by Jutila--Motohashi [8].

Note that (1.7) and (1.9) do not seem to apply one another,
and that for the
derivation of (1.8) from (1.9) the non-negativity of $H_j(\hf)$
is not needed.

\medskip
Our new results on sums of sums of $\a_jH_j^3(\hf)$ and $\a_jH_j^4(\hf)$
are contained in

\bigskip
THEOREM 1. {\it We have}
$$
\sum_{\k_j\le K} \a_jH_j^3(\hf) = K^2P_3(\log K) + O(K^{5/4}\log^{37/4}K),
\leqno(1.10)
$$
{\it where $P_3(x)$ is a polynomial of degree three with leading
coefficient $4/(3\pi^2)$, whose remaining
coefficients may be explicitly evaluated.}

\bigskip
THEOREM 2. {\it We have}
$$
\sum_{\k_j\le K} \a_j\H = K^2P_6(\log K) + O(K^{3/2}\log^{25/2}K),
\leqno(1.11)
$$
{\it where $P_6(x)$ is a polynomial of degree six with leading
coefficient $16/(15\pi^4)$, whose  remaining
coefficients may be explicitly evaluated.}

\bigskip
The proofs of (1.10) and (1.11), which will be given in subsequent
sections, depend on several ingredients. Besides the transformation
formulas for sums of $\a_jH_j^k(\hf)$ (see Section 3), two  salient ones
are the short interval bounds (1.7) and (1.9),
and the estimates for the sixth
and eighth moments of $|\zt|$. Indeed, it is a deep and beautiful
fact that sums of $\a_jH^k_j(\hf)$ and moments of $|\zt|^{2k}\;(k\in\NN)$
are closely related, at least for $k \le 4$. Both quantities tend
to increase in complexity as $k$ increases. One of the reasons why
Motohashi was able to get the sharp error term $O(T\log^6T)$ in (1.2)
was that the continuous part of his relevant formula, namely the integral
on the left-hand side of (1.12) below, contained $|\zt|^4$.
However, for $\int_0^T|\zt|^4\d t$ we know that the correct order
of magnitude is $T\log^4T$, and actually the asymptotic formula with
error term $O(T^{2/3}\log^CT)$ is known (see e.g., [4] and [15]).
Unfortunately, to this day such type of result is not known for any
power moment of $|\zt|$ greater than the fourth.

As to the true order of sums of
$\a_j H^k_j(\hf)$, perhaps it is true that, for $k\in \NN$ fixed,
$$
\sum_{\k_j \le T}\a_j H_j^k(\hf)
+ {2\over{\pi}}\int_{0}^T {|\zt|^{2k}\over|\z(1+2it)|^2}\d t
= T^2P_{{1\over2}(k^2-k)}(\log T)
+ O(T^{1+c_k+\e}),\leqno(1.12)
$$
where $P_{{1\over2}(k^2-k)}(z)$ is a suitable polynomial of degree
${1\over2}(k^2-k)$ in $z$ whose coefficients depend on $k$,
and $0 \le c_k < 1$; perhaps even $c_k = 0$ is true. We actually
have $c_2 = 0$ in view of (1.2), and from the proofs of Theorem 1 and
Theorem  2 it follows that we may take $c_3 = 1/7, c_4 = 1/3$.
For example, (5.5) and (5.6) (for $k=4$)
clearly show why the left-hand side of (1.12) appears, and in view of
$H_j(\hf) \ge 0$ it is positive. It would be interesting to evaluate
(or estimate) the sum in (1.12) when $k=1$ and $k\ge5$. The case $k=1$
will be briefly discussed at the end of the paper, while $k\ge5$ lies outside
the scope of this work. However the latter case is of potential
importance since it could yield upper bounds for the $2k$-th
moment of $|\zt|$. Namely if for some $k\ge 6$ the right-hand side
of (1.12) is bounded by $T^{2+\e}$, this would essentially give a
bound at least as strong as the (known) twelfth moment of $|\zt|$
(see (4.2)). If this bound holds for every $k$, then this implies both
$H_j(\hf) \ll_\e \k_j^\e$ and the Lindel\"of hypothesis that $\zt
\ll_\e |t|^\e$. It is yet unknown what is the connection between these
two conjectures, namely whether one of them implies the other one.

\medskip
Conjectures for moments of various $L$-functions have
been recently proposed by considerations from Random matrix
theory (see J.B. Conrey [1] and the comprehensive work by J.B. Conrey, D.W.
Farmer, J.P. Keating, M.O. Rubinstein and N.C. Snaith [2]). In all cases
which can be predicted by this theory and where the asymptotic formula
in question was rigorously proved, the main terms coincide. In our
context this theory says that one should have
$$
\sum_{\k_j \le K}\a_j H_j^k(\hf) = K^2P_{{1\over2}(k^2-k)}(\log K) + o(K^2)
\leqno(1.13)
$$
for $k\in\NN$ fixed. The leading coefficient of
$P_{{1\over2}(k^2-k)}(x)$ equals
$$
d_k \;=\; {a_k g_k\over\pi^2\left({k(k-1)\over2}\right)!}.\leqno(1.14)
$$
In the notation of this theory $g_k$ is the so-called {\it geometric}
part. In our case it is
$$
g_k = \left(\hf k(k-1)\right)!\,2^{k(k+1)/2-1}\prod_{j=1}^{k-1}
{j!\over(2j)!},  \leqno(1.15)
$$
so that $g_1 = 1, g_2 = 2, g_3 = 8, g_4 = 128$. The constant $a_k$ is
the {\it arithmetic} part. It equals
$$
a_k = \prod_p\left(1-{1\over p}\right)^{k(k-1)/2}
\sum_{j=0}^\infty {k+j-1\choose j} {k+j-2\choose j}{1\over(j+1)p^j}.
\leqno(1.16)
$$
We have $a_1 = a_2 = a_3 = 1, a_4 = 1/\z(2) = {6\over\pi^2}$. In general,
$a_k$ can be expressed in terms of hypergeometric functions.
Note that
$$
\sum_{j=0}^\infty  {k+j-1\choose j} {k+j-2\choose j}{x^j\over j+1}
\qquad(|x| < 1)
$$
is a rational function of $x$ whose denominator is $(1-x)^{2k-3}$ and
numerator is 1 for $k =2,\,3$, and is equal to $1+x \,(k=4),\, 1+3x+x^2\,
(k=5)$ etc. This shows that, for $k\ge5$, $a_k$ will not be expressible in a
simple closed form, but as an Euler product over the primes.
We have the values $d_1 = 1/\pi^2, d_2 = 2/\pi^2, d_3 = 4/(3\pi^2),
d_4 = 16/(15\pi^4)$, which coincide for $k = 2,3,4$ with the ones that
follow from (1.2), Theorem 1 and Theorem 2.  Note that Random matrix theory
also predicts the asymptotic formula for the sum in (1.13) without the
normalizing factor $\alpha_j$. The shape of the conjectured formula
will be similar to the above one, only the constants will be different,
and somewhat more complicated. Unfortunately, the methods at hand permit
one to deal only with the sum in (1.13).
\medskip
{\bf Acknowledgement}. I wish to thank Prof. Brian Conrey and
Prof. Matti Jutila for valuable remarks.
\bigskip
\head
2. Kuznetsov's work on sums of $\H$
\endhead

N.V. Kuznetsov's preprint [12] states as the main result (Theorem 1 on p. 5)
the asymptotic formula
$$
\sum_{\k_j\le T}\a_j\H = T^2P_6(\log T) + O(T^{4/3+\e}),\leqno(2.1)
$$
where $P_6(x)$ is a polynomial in $x$ of degree six whose leading
coefficient is equal to $2^9/(15\pi^3)$. This is actually
stronger than our (1.11) of Theorem 2. Unfortunately, Kuznetsov
did not prove (2.1), and even the leading coefficient of $P_6(x)$
is not correctly stated (it equals $16/(15\pi^4)$, see Section
9 for details). We shall analyze his preprint and substantiate
our claim, using certain valid parts of his work, namely the derivation
of the main term to shorten the proof of our Theorem 2.
A complete list of misprints, errors etc. of [12] is not given,
but just some of the important ones will be stated here. Further
discussion concerning [12] will be given in subsequent sections.

\smallskip
Page  8, line after (21) it is not shown  why $\hat{\psi}(2w)$ is regular for
$\R w > -5/2$, which is claimed in the text. Namely
$$
\hat{\psi}(2w) = \int_{-\infty}^\infty
\pi^{-1}2^{2w-1}\G(w-iu)\G(w+iu)
h(u)u\,{\sinh}\,(\pi u)\d u\;(\R w > 0),\leqno(2.2)
$$
where ($Q \asymp T^{1/3}$)
$$
\eqalign{
h(r) &= q(r)\left\{\exp\left(-\left({r-T\over Q}\right)^2\right)
+ \exp\left(-\left({r+T\over Q}\right)^2\right)\right\},\cr
q(r) &= {(r^2 + {\txt{1\over4}})(r^2 + {\txt{9\over4}})\over
(r^2 + {\txt{1\over4}})(r^2 + {\txt{9\over4}}) + 626},\cr}\leqno(2.3)
$$
so that $h(r)$ is even, regular for $|\I w| \le 3$,
$h(\pm {i\over2}) = h(\pm {3i\over2}) = 0$,
and $h(r)$ decays like $\exp(-c|r|^2)$. To analyze the function
$\hat{\psi}(2w)$, note that from
$$
\G(z)\G(1-z) = {\pi\over\sin(\pi z)}\leqno(2.4)
$$
one obtains the identity
$$
\G(w+iu)\G(w-iu) = {\pi i\over2\sinh(\pi u)\cos(\pi w)}
\left\{{\G(w+iu)\over\G(1-w+ iu)} - {\G(w-iu)\over\G(1-w - iu)}\right\}.
$$
Since $h(r)$ is even, this gives
$$
\hat{\psi}(2w) \;=\; {i2^{2w}\over\cos(\pi w)}h^*(w),\leqno(2.5)
$$
where (see Y. Motohashi [14, eq. (2.12)])
$$
h^*(s) \;:=\; \int_{-\infty}^\infty uh(u){\G(s+iu)\over\G(1-s+iu)}\d u
= -\int_{(0)}wh(iw){\G(s+w)\over\G(1-s+w)}\d w\leqno(2.6)
$$
is regular for $\R s > 0$, where $\int_{(\a)} $ means integration over
the line $\R w = \a$. If $h(r)$ is entire (cf. (1.6)), then in (2.6)
the line of integration may be shifted to  $\R w = C > 0$. Thus
$h^*(s)$  is seen to be regular for $\R s > -C$, and since $C$
may be  arbitrary, it follows
that $h^*(s)$ is entire and of polynomial growth in $|s|$ for $\s$
in a fixed strip. In the case of (2.3) $h^*(s)$ is regular at least
for $\R s > -3$, and we have
$h^*(\pm\hf) = h^*(\pm{3\over2}) = 0$. For example, by using taking
$\R w = 2$ in (2.6) and using
the functional equation $s\G(s) = \G(s+1)$ one obtains
$$
{\G(-{3\over2}+w)\over\G({5\over2}+w)} = {1\over(w^2-{1\over4})
(w^2 - {9\over4})}.
$$
Thus this cancels with the corresponding factor of $h(iw)$, and
$h^*(-{3\over2}) = 0$ follows since $h(r)$ is even. Likewise it follows
that $h^*(n+\hf) = 0\;(n\in\NN)$, hence $\hat{\psi}(2w)$ is indeed
regular for $\R w > -5/2$, the first pole at $w = -5/2$ coming from
the zero of $\cos(\pi w)$ in the numerator in (2.5).

\smallskip
Page 9, in the formulation of Theorem 2 the numbering (27) is missing,
and the condition (contradicting $\R  \mu = \R \nu = \hf$)
$\R \mu, \R \nu \not = \hf$ should be $\mu, \nu \not = \hf$.

\smallskip
More importantly, Kuznetsov did not prove Theorem 2 (which yields
the spectral decomposition for the sum in (3.2), and is the basis
of [12]) in [10] as he claimed.
The result was used there in his unsuccessful attempt to prove
the eighth moment for the Riemann zeta-function, namely
$$
\int_0^T |\zt|^8\d t \ll T\log^CT.
$$
The same formula was also used in [11] in his failure to prove the
Lindel\"of hypothesis that $\zt \ll_\e |t|^\e$. A corrected version
of the formula is due to Y. Motohashi [13] in 1991, and
recently this was updated and improved in [16]. Hence due
to Motohashi's work [16] this important obstacle in
dealing with the asymptotic evaluation of the sum in
(2.1) has been removed but, unfortunately, this is
not the only shortcoming of [12] as will be clear from the sequel.

\smallskip
Page 10, l. 8. Kuznetsov chooses $s = \nu = \rho = \mu = \hf$,
which violates the assumptions of Theorem 2, without mentioning
that first one has to take $\mu = \hf + it, \nu = \hf + i\tau$
and then to take $t, \tau\to 0$. In (30), in the first line,
$1-2s$ should be $1-2\rho$.

\smallskip
Page 14, in l. 10 (32) should be (38), in (44) $4^5 = 1012$ is false.

\smallskip
Page 24, l. -5,6  it should be sh$\,\eta/2 = {r\over T}$.

\smallskip
Page 25. l. 2,4 of (91), $\xi$ is repeatedly written in place of $\z$.
Formula (92) is incorrect, detailed discussion will be given below
in Section 4.
In (93), on the right-hand side, $Q$ is missing twice. In (94), in
the exponent in the $O$-term, $ve$ should be replaced by $\e$.
In (95), $dt$ should be $dT$, $r^{3/2}$ should be $\k_j^{3/2}$.
Line below (91), $\xi$ should be (6).

\smallskip
Page 26, in (97) $Q$ is missing once on both sides, $ve$ should be $\e$.

\bigskip
\head
3. Formulas for products of three and four Hecke series
\endhead

The essence to the approach of dealing with sums of $H_j^3(\hf)$
and $\H$ are the transformation formulas for the sums
$$
{\Cal C}(K,G) := \sum_{j=1}^\infty \a_jH_j^3(\hf)h_0(\k_j)
\leqno(3.1)
$$
with    $h_0(r) $ given by (1.6), and
$$
\sum_{j=1}^\infty \a_j\H h(\k_j),
\leqno(3.2)
$$
with $h(r)$ given by (2.3). The notation in (3.2) corresponds to
Motohashi [15], while that of (2.3) is from Kuznetsov  [12]. We shall
adhere to this for practical reasons, but of course it would have
been possible to use $h_0(\k_j)$ instead of $h(\k_j)$ etc. To obtain
transformation formulas for the weighted sums (this facilitates
the resolution of the problems involving analytic continuation) one
starts from general expressions, namely $H_j(u)H_j(v)H_j(\hf)$ in
(3.1) and $H_j(u)H_j(v)H_j(w)H_j(z)$ in (3.2) in the region of
absolute convergence. In the former one
replaces $H_j(\hf)$ by an approximate
functional equation (e.g., [15, Lemma 3.9]) which reduces it to
suitable sums of $t_j(f)f^{-1/2}$. The product of two Hecke
series is transformed by the use of the identity (in the region of
absolute convergence; see [15, (3.2.7)])
$$
H_j(s)H_j(s-\a) = \z(2s-\a)\sum_{n=1}^\infty \s_\a(n)t_j(n)n^{-s}
\quad(\s_\a(n) = \sum_{d|n}d^\a),
$$
which is the analytic equivalent of the multiplicativity of the
arithmetic function $t_j(n)$, namely (see e.g., [15, eq. (3.1.14)])
$$
t_j(m)t_j(n) \= \sum_{d|(m,n)}t_j\left({mn\over d^2}\right). \leqno(3.3)
$$
 After this there is
summation of $t_j(m)t_j(f)$
in both cases, which is effected by applying the Kuznetsov trace
formula (see [15, Theorem 2.4]). It is here that delicate questions
of analytic continuation arise. In [7] M. Jutila used a variation
of this approach in proving (1.9). Namely he used ([14, pp. 266-267]
or [15, Lemma 3.8]) Motohashi's formula for
$$
\sum_{j=1}^\infty \a_jH^2_j(\hf)t_j(f)h(\k_j)\qquad(f \in \NN),
$$
combined with his explicit expression for $H_j^2(\hf)$ (see (9.2)).

\smallskip
We shall present now Motohashi's explicit formula
for sums of $H_j^3(\hf)$, needed for the proof of Theorem 1.
We have (see (3.1)) with
$\lambda = C\log K\;(C > 0)$  (this is [15, (3.5.18)], with the
extraneous factor $(1 - (\k_j/K)^2)^\nu$ omitted)
$$\eqalign{
{\Cal C}(K,G) &= \sum_{f\le3K}f^{-{1\over2}}\exp\Bigl(-{\bigl({f\over K}
\bigr)}^\lambda\Bigr){\Cal H}(f;h_0)\cr&
- \sum_{\nu=0}^{N_1}\,\sum_{f\le3K}f^{-{1\over2}}U_\nu(fK)
{\Cal H}(f;h_\nu) + O(1),\cr}                    \leqno(3.4)
$$
with ($h_0(r)$ is given by (1.6))
$$
h_\nu(r) \= h_0(r)\left(1 - \left({r\over K}\right)^2\right)^\nu
\quad(\nu = 0,1,2\ldots),\leqno(3.5)
$$
$$
{\Cal H}(f;h) \= \sum_{\nu=1}^7{\Cal H}_\nu(f;h),
$$
$$
{\Cal H}_1(f;h) \= -2\pi^{-3}i\left\{(\gamma - \log(2\pi{\sqrt f}))
({\hat h})'(\hf) + {\txt{1\over4}}({\hat h})''(\hf)\right\}d(f)
f^{-{1\over2}},
$$
$$
{\Cal H}_2(f;h) \= \pi^{-3}\sum_{m=1}^\infty m^{-{1\over2}}d(m)d(m+f)
\Psi^+({m\over f};h)\quad\Big(d(n) = \sum_{\delta|n}1\Big),
$$
$$
{\Cal H}_3(f;h) \= \pi^{-3}\sum_{m=1}^\infty (m+f)^{-{1\over2}}
d(m)d(m+f)\Psi^-(1+ {m\over f};h),                      \leqno(3.6)
$$
$$
{\Cal H}_4(f;h) \= \pi^{-3}\sum_{m=1}^{f-1}m^{-{1\over2}}d(m)d(f-m)
\Psi^-({m\over f};h),
$$
$$
{\Cal H}_5(f;h) \= -(2\pi^3)^{-1}f^{-{1\over2}}d(f)\Psi^-(1;h),
$$
$$
{\Cal H}_6(f;h) \= -12\pi^{-2}i\s_{-1}(f)f^{1\over2}h'(-\hf i),
$$
$$
{\Cal H}_7(f;h) \= -\pi^{-1}\int_{-\infty}^\infty{|\z(\hf+ir)|^4\over
|\z(1+2ir)|^2}\s_{2ir}(f)f^{-ir}h(r)\d r,
$$
where (see (2.6))
$$
\Psi^+(x;h) \= \int_{(\b)}\G^2(\hf-s)\tan(\pi s)h^*(s)x^s\d s,
$$
and
$$
\Psi^-(x;h) \= \int_{(\b)}\G^2(\hf-s){h^*(s)\over
\cos(\pi s)}x^s\d s
$$
with $-{3\over2} < \b < {\hf}$. In (3.4) $N_1$ is a sufficiently
large integer and
$$
U_\nu(x) = {1\over2\pi i\lambda}\int_{(-\lambda^{-1})}
(4\pi^2K^{-2}x)^wu_\nu(w)
\G({w\over\lambda})\d w \ll \left({x\over K^2}\right)^{-{C\over\log K}}
\log^2K,\leqno(3.7)
$$
where $u_\nu(w)$ is a polynomial in $w$ of degree $\le 2N_1$, whose
coefficients are bounded. A prominent feature of Motohashi's explicit
expression for ${\Cal C}(K,G)$ is that it contains series and integrals
with the classical divisor function $d(n)$ only, with no quantities
from spectral theory. Therefore the problem of evaluating
 ${\Cal C}(K,G)$ is a problem of classical analytic
number theory.

\medskip
As for (3.2), we adopt the notation of [12], primarily since we intend
to correct Kuznetsov's proof. As already stated, a correct and
rigorous proof of the spectral decomposition for (3.2) is given by
Y. Motohashi [13] and [16]. The formulation is technically complicated,
and for the sake of brevity will not be reproduced here.

\bigskip
\head
4. The asymptotic formula for sums of $\H$
\endhead

We shall provide in fact two completely different proofs
of Theorem 2. The first is obtained by correcting and simplifying
the proof given by N.V. Kuznetsov in [12].
The second approach consists of elaborating the method
of M. Jutila [7], used in the proof of the bound (1.9), which is
one of the crucial ingredients in the proof   of Theorem 2.
It will be outlined in Section 9.

We shall begin  now with the proof of (1.11) of Theorem 2,
correcting and simplifying [12].
We remark first that one obtains (1.11) from
$$
\sum_{\k_j\le T}\a_j\H + O\left(\log^2T\int_0^T|\zt|^8\d t \right)
= T^2P_6(\log T) + O(T^{{4/3}+\e}).\leqno(4.1)
$$
Namely one has (e.g., see [3]) the bounds
$$
\int_0^T|\zt|^4\d t \ll T\log^4T,\qquad
\int_0^T|\zt|^{12}\d t \ll T^2\log^{17}T.\leqno(4.2)
$$
Hence by the Cauchy-Schwarz inequality for integrals it
follows from (4.2) that
$$
\int_0^T|\zt|^8\d t \ll T^{3/2}\log^{21/2}T, \leqno(4.3)
$$
which is still the sharpest known upper bound estimate for
the integral in (4.3). In [12] N.V. Kuznetsov assumed that the bound
$$
\int_0^T|\zt|^8\d t \ll T\log^{C}T \leqno(4.4)
$$
holds for some $C>0$. This is what he claimed to have proved in [10].
Although he never officially withdrew the claim (the proof was faulty),
this fact was mentioned in the review in the Zentralblatt (Zbl.745.11040).
The asymptotic formula (4.1) shows clearly that one cannot attain the
exponent $4/3+\e$ in (4.1) unless it is attained in (4.3). This, however,
would be a big achievement in zeta-function theory.

The plan of the
proof is as follows: from the fundamental formula for sums of
products of four of Hecke series ([12, Theorem 2] or [16, Theorem])
one obtains first the  formula
$$
\eqalign{&
\sum_{j\ge1}\a_j\H h(\k_j)
+ {2\over\pi}\int\limits_0^\infty{|\zeta(\hf+ir)|^8\over|\z(1+2ir)|^2}
h(r)\d r\cr&
+ \sum_{k\ge 12,k\equiv0({\roman{mod}})2}g(k)\sum_{1\le j\le\nu_k}
\a_{j,k}H_{j,k}^4(\hf)\cr&
= \sum_{j\ge1}\a_j\H(h_0(\k_j) + \e_jh_1(\k_j))
+ {1\over\pi}\int\limits_{-\infty}^\infty
{|\zeta(\hf+ir)|^8\over|\z(1+2ir)|^2} (h_0(r)+ h_1(r))\d r\cr&
+ R + O(Q\log^6T).
\cr}\leqno(4.5)
$$
Here $h(r)$ is given by (2.3), the quantities in
$$
\sum_{k\ge 12,k\equiv0({\roman{mod}})2}g(k)\sum_{1\le j\le\nu_k}
\a_{j,k}H_{j,k}^4(\hf),\leqno(4.6)
$$
which are associated with holomorphic cusp forms are precisely
defined in [12] or [15],
$$
Q = T^{1/3}, \leqno(4.7)
$$
$$\eqalign{&
g(k) = {1\over2\pi^3 i}\int_{(\delta)}{\G(k-\hf+w)\over\G(k+\hf +w)}
\G^4(\hf - w)\sin(\pi w)\hat{\psi}(2w)\d w,
\cr&
h_0(r) = {1\over\pi^3 i}\int_{(\delta)}\G(w+ir)\G(w-ir)
\G^4(\hf - w)\sin(2\pi w)\hat{\psi}(2w)\d w,
\cr&
h_1(r) = {1\over2\pi^3 i}\int\limits_{(\delta)}\G(w+ir)\G(w-ir)
\G^4(\hf - w)\cosh(\pi r){\sin^2(\pi w)+1\over\cos(\pi w)}
\hat{\psi}(2w)\d w,\cr}
\leqno(4.8)
$$
where $\hat\psi$ is given by (2.2)  and $\delta>0$ is a small constant.
The choice of $Q$ in (4.7) seems optimal, and any improvements
(namely $Q = T^\a$ with $\a < 1/3$) will require the use of new methods.
Actually, instead of (4.7) the correct choice of $Q$ is $Q = CT^{1/3}$
with some $C > 0$, since we shall integrate (4.5) over the interval
$[T_0,\,2T_0]$, so $Q$ should ultimately depend on $T_0$ and not on $T$;
e.g., one can take $Q = T_0^{1/3}$ (this fact is not mentioned in [12]).
The symbol $R$ in (4.5) stands for the residual (main) terms. This has
been calculated by Kuznetsov in [12] to  be equal to
$$
\sum_{k=0}^6a_k{\hat{\psi}}^{(6-k)}(1)  + O(Q\log^6T)
\qquad\bigl(a_0 = {1024\over15\pi^3}\bigr).\leqno(4.9)
$$
It can be shown that the contribution  of (4.6) is $O(Q\log^6T)$
(note that the sum with $\a_{j,k}H_{j,k}^4(\hf)$ is easily majorized; see
[15]) and so is the contribution of $g(k)$ and $h_0$ in (4.5)
(see (6.2) and (6.3)). What remains then is the basic formula
$$
\eqalign{&
\sum_{j\ge1}\a_j\H\exp\left(-\left({\k_j-T\over Q}\right)^2\right)
+ {2\over\pi}\int\limits_0^\infty{|\zeta(\hf+ir)|^8\over|\z(1+2ir)|^2}
\exp\left(-\left({r-T\over Q}\right)^2\right)\d r\cr&
= \sum_{k=0}^6a_k{\hat{\psi}}^{(6-k)}(1) + \sum_{j\ge1}\a_j\H{\tilde h}(\k_j)
\cr& +\,
{1\over\pi}\int\limits_0^\infty
{|\zeta(\hf+ir)|^8\over|\z(1+2ir)|^2}{\tilde h}(r)\d r
+ O(Q\log^6T),
\cr}\leqno(4.10)
$$
where  ${\tilde h}(r)$ is the oscillatory integral transform
obtained by replacing $\sin^2(\pi w) + 1$ in the definition of $h_1(r)$
(see (4.8)) by $\sin^2(\pi w)$. The terms containing this function  will
be small, while $\hat \psi$
will give rise to the main term $T^2P_6(\log T)$ in (4.1).

In the relevant range one has (this follows from Kuznetsov's
Lemma 4.7)
$$
\tilde{h}(r) \ll Qr^{-1/2}\exp(-CQ^2r^2T^{-2})\qquad(C>0). \leqno(4.11)
$$
Hence by the non-negativity of the integral on the left-hand side
of (4.10), (1.3) and (4.11) it follows that
$$
\eqalign{&
\sum_{j\ge1}\a_j\H\exp\left(-\left({\k_j-T\over Q}\right)^2\right) \cr&
\ll QT\log^6T + Q\sum_{\k_j\le TQ^{-1}\log T}\a_j\H \k_j^{-1/2}\cr&
\ll QT\log^6T + T^{3/2}Q^{-1/2}\log^{16}T
\ll T^{4/3}\log^{16}T.\cr}\leqno(4.12)
$$
Observe that (4.12) is a sharpened variant, in view of (4.7), of
(1.9), as it gives (1.9) with the right-hand side replaced by
$K^{4/3}\log^{16}K$.
By using (4.3) and (4.11) it follows that the integral  on the
right-hand side of (4.5) is $\ll T^{1+\e}$. Also note that we have
($Q = T^{1/3}$)
$$
\int_0^\infty{|\zeta(\hf+ir)|^8\over|\z(1+2ir)|^2}
\exp\left(-\left({r-T\over Q}\right)^2\right)\d r \ll T^{4/3}\log^{16}T,
\leqno(4.13)
$$
which can be easily obtained from the mean square bounds for
$\zt$ over short intervals (see [3, Chapter 15]) and the classical
bound $\zt \ll |t|^{1/6}$. One also has to use e.g. the standard bound
$$
{1\over|\z(1+it)|} \;\ll\; \log |t|.\leqno(4.14)
$$

\smallskip
After these considerations it remains to integrate the
basic formula (4.10) over $T$
from $T_0$ to $2T_0$ and then to replace $T_0$ by $T_02^{-j}$, and
sum the resulting expressions for $j\ge1$. This will lead
to (4.1). The technical details are given in the next section,
as well as the calculation of the main term.

We have restrained
ourselves from analyzing the difficult lemmas of [12, Section 4],
especially of the Lemma 4.7 which claims an asymptotic formula
for the crucial function ${\tilde h}(r)$ appearing in (4.10).
The function $h_1(r)$ in [12, Lemma 3.2] is first transformed
into a complicated expression involving  the hypergeometric function.
This is said to follow from the use of Parseval's formula for Mellin
transforms. The author was unable to follow the proof of
Lemma 4.7, which claims an asymptotic expansion of ${\tilde h}(r)$.
However, this  asymptotic expansion will be proved,
in Section 6,   by a method which is
different and simpler than Kuznetsov's.

\bigskip
\head
5. Integration of the basic formula and the main term
\endhead

We shall deal first with the main term in (1.11). One way
to obtain this expression is to go through Kuznetsov's
paper [12]. Therein he claimed (eq. (92) on p. 25) that
$$
\hat{\psi}^{(m)}(1) = {2\over\sqrt{\pi}}QT\log^mT
\cdot\left(1 + O\left({1\over T}\right)\right)
\qquad(m = 0,1,2,\cdots),\leqno(5.1)
$$
where $\hat{\psi}$ is defined by (2.2). On the right-hand side of (4.10)
there appears
$$
\sum_{k=0}^6a_k{\hat{\psi}}^{(6-k)}(1),
$$
which will give rise to the main term $K^2P_6(\log K)$ in (2.1).
Hence we have to evaluate explicitly $\hat{\psi}^{(m)}(1)$ for
$m = 0,\cdots,6$.

\smallskip
{\it The case $m=0$.} From (2.2) we have, on using (2.4)
and recalling that $h(r)$ is given by (2.3),
$$
\eqalign{
\hat{\psi}(1) & = {2\over\pi}\int_0^\infty \G(\hf+iu)\G(\hf-iu)uh(u)
\,{\sinh}\,(\pi u)\d u\cr&
= 2\int_0^\infty\,{uh(u)\over\sin\pi(\hf+iu)}\,{\sinh}\,(\pi u)\d u
= 2\int_0^\infty\,uh(u)\,{\tanh}\,(\pi u)\d u\cr&
= 2\int_{T-Q\log T}^{T+Q\log T} \exp(-(u-T)^2Q^{-2})u\bar{h}(u)
\,{\tanh}\,(\pi u)\d u + O({\roman e}^{-\hf\log^2T}),\cr}
$$
where
$$
\bar{h}(r) \= 1 + O(r^{-4}).
$$
Change of variable $u = T + Qx$ gives then
$$
\eqalign{
\hat{\psi}(1) &= 2Q\int_{-\log T}^{\log T}{\roman e}^{-x^2}(T+Qx)
\,{\tanh}\,\pi (T + Qx)\d x  + O(1)\cr&
= 2Q(\sqrt{\pi}T(1 + O(1/T)) + O(Q)) = 2\sqrt{\pi}QT\left(1 +
O\left({Q\over T}\right)\right).\cr}\leqno(5.2)
$$

\smallskip
{\it The case $m \ge 1$}. We need the formula (see e.g., [4, p. 272])
$$
{\G^{(k)}(s)\over\G(s)} = \sum_{j=0}^k b_{j,k}(s)\log^js
+ c_{-1,k}s^{-1} +\cdots + c_{-1,r}s^{-r} + O_r(|s|^{-r-1})\leqno(5.3)
$$
for fixed integers $k \ge 1,\,r \ge 0$, where each of the functions
$b_{j,k}(s) \;(\sim b_{j,k}$ for a suitable constant  $b_{j,k}$ as
$s\to\infty$) has an asymptotic expansion in non-positive powers of
$s$. As in the case $m=0$ the main contribution to $\hat{\psi}^{(m)}(2w)$
will come from an interval of length $\ll Q\log T$, when $w$ lies
in a neighbourhood of $\hf$. Namely we have
$$\eqalign{
2^m\hat{\psi}^{(m)}(2w) &= {1\over\pi}\int\limits_{T-Q\log T}^{T+Q\log T}
{\d^m\over\d w^m}\left(2^{2w}\G(w+iu)\G(w-iu)\right)uh(u){\sinh}\,(\pi u)
\d u \cr& + O({\roman e}^{-\hf\log^2T}).\cr}
$$
To calculate the derivatives in the above integral we apply Leibniz's rule.
We have to evaluate ($r = 0 ,1,\cdots ,m)$
$$
{\d^r\over\d w^r}\G(w+iu)\G(w-iu)\Biggl|_{w=\hf},\qquad u = T + O(Q\log T).
$$
By using (5.2), (5.3)  and (2.4) it is seen that this expression equals
$$
\eqalign{&
\sum_{j=0}^r{r\choose j}\G^{(j)}(\hf+iu)\G^{(r-j)}(\hf-iu)\cr&
={\pi\over{\cosh\,}(\pi u)}\left(\sum_{\ell=0}^r d_{\ell,r}\log^ru
+ O(T^{-1}\log^rT)\right)
\cr}
$$
with suitable constants $d_{\ell,r}$. Proceeding as in the case
$m=0$, we obtain
$$
\hat{\psi}^{(m)}(1) = QT\left(\sum_{j=0}^m c_{j,m}\log^mT
+ O_m(QT^{-1}\log^mT)\right)\quad(m\in\NN)\leqno(5.4)
$$
with suitable constants $c_{j,m}$, which may be explicitly evaluated
($c_{m,m} = 2^{2-m}\sqrt{\pi}$). From (5.2) and (5.4) we see
that Kuznetsov's claim (5.1) is incorrect.

\medskip
Now we integrate (4.10) over $T$ from $T_0$ to $2T_0$, taking
$Q = T_0^{1/3}$ (cf. (3.7)), which clearly may be done. We have first
$$
\eqalign{&
\sum_{j\ge1}\a_j\H\int_{T_0}^{2T_0}
\exp\left(-\Bigl({\k_j-T\over Q}\Bigr)^2\right)\d T\cr&
= \sum_{T_0-Q\log T_0\le\k_j\le2T_0+Q\log T_0}\a_j\H\int_{T_0}^{2T_0}
\exp\left(-\Bigl({\k_j-T\over Q}\Bigr)^2\right)\d T + o(1).
\cr}
$$
By change of variable and  (4.12) (or (1.9)) the sum on the right-hand
side equals
$$
\eqalign{&
Q\sum_{T_0-Q\log T_0\le\k_j\le2T_0+Q\log T_0}\a_j\H
\int_{(T_0-\k_j)/Q}^{(2T_0-\k_j)/Q}{\roman e}^{-x^2}\d x\cr&
= O(Q^2T_0^{1+\e}) + Q\sum_{T_0+Q\log T_0\le\k_j\le2T_0-Q\log T_0}\a_j\H
\int_{(T_0-\k_j)/Q}^{(2T_0-\k_j)/Q}{\roman e}^{-x^2}\d x\cr&
= O(Q^2T_0^{1+\e}) + \sqrt{\pi}Q
\sum_{T_0+Q\log T_0\le\k_j\le2T_0-Q\log T_0}\a_j\H  +
O({\roman e}^{-\hf\log^2T})\cr&
= \sqrt{\pi}Q\sum_{T_0\le\k_j\le2T_0}\a_j\H + O(Q^2T_0^{1+\e}).
\cr}\leqno(5.5)
$$
In a similar fashion, by using (4.3) and (4.13), it follows that
$$
\eqalign{&
\int_{T_0}^{2T_0}\int_0^\infty{|\zeta(\hf+ir)|^8\over|\z(1+2ir)|^2}
\exp\left(-\left({r-T\over Q}\right)^2\right)\d r\,\d T\cr&
=\sqrt{\pi}Q\int_{T_0}^{2T_0}{|\zeta(\hf+ir)|^8\over|\z(1+2ir)|^2}\d r
+ O(Q^2T^{1+\e})\cr&
\ll Q\log^2T_0\int_{T_0}^{2T_0}|\zeta(\hf+ir)|^8\d r + Q^2T^{1+\e}
\ll QT_0^{3/2}\log^{25/2}T_0.   \cr}\leqno(5.6)
$$
To bound the second sum on the right-hand side of (4.10)
we use Lemma 4.7 of [12], or the discussion on $h_1(r)$ in Section 6.
We need especially the terms  $(T-r)\log(T-r) - (T+r)\log(T+r)$ in (6.10),
in conjunction with the first derivative test (Lemma 2.1 of [3])
and (1.3). The derivative in question is $\gg r/T$, and we shall obtain
(${\tilde h}(\k_j) = {\tilde h}(\k_j,T)$)
$$
\eqalign{&
\sum_{j\ge1}\a_j\H\int_{T_0}^{2T_0}\tilde{h}(\k_j,T)\d T\cr&
\ll \sum_{j\ge1}\a_j\H QT_0\k_j^{-3/2}\exp\left(-{Q^2\k_j^2\over4T_0^2}
\right) + 1\cr&
\ll QT_0\sum_{\k_j\le T_0Q^{-1}\log T}\a_j\H \k_j^{-3/2} + 1 \cr&
\ll QT_0^{{4\over3}+\e}.
\cr}\leqno(5.7)
$$
Finally from (5.2) and (5.4) we have
$$\eqalign{&
\int_{T_0}^{2T_0}\sum_{k=0}^6a_k{\hat{\psi}}^{(6-k)}(1) \d T
= Q\int_{T_0}^{2T_0}T\sum_{k=0}^6a_k\sum_{j=0}^{6-k} e_{j,m}\log^mT\d T +
 O(QT_0^{{4\over3}+\e})\cr& =
QT^2\sum_{k=0}^6 f_{k}\log^kT\Bigg|_{T_0}^{2T_0} +O(QT_0^{{4\over3}+\e}).
\cr}
\leqno(5.8)
$$
with effectively computable constants $e_{j,m}$ and
$f_{k}$. Therefore (4.1) will follow
from (5.5)--(5.8) when we divide by $Q$, replace $T_0$ by $T_02^{-j}$
and sum over $j$.

\bigskip
\head
6. The estimates for the oscillatory terms
\endhead
In this section we shall complete the proof of Theorem 2 by
estimating the oscillatory functions defined by (4.8). We shall
use the function $h^*(s)$, defined by (2.5)--(2.6) to simplify the
functions in (4.8). We obtain
$$\eqalign{&
g(k) = {1\over\pi^3}\int_{(\delta)}2^{2w-1}{\G(k-\hf+w)\over\G(k+\hf +w)}
\G^4(\hf - w)\tan(\pi w)h^*(w)\d w\quad(k \ge 12),
\cr&
h_0(r) = {1\over\pi^3}\int_{(\delta)}2^{2w+1}\G(w+ir)\G(w-ir)
\G^4(\hf - w)\sin(\pi w)h^*(w)\d w,
\cr&
h_1(r) = {1\over2\pi^3}\int\limits_{(\delta)}2^{2w}\G(w+ir)\G(w-ir)
\G^4(\hf - w)\cosh(\pi r){\sin^2(\pi w)+1\over\cos(\pi w)}h^*(w)\d w,\cr}
\leqno(6.1)
$$
where $\delta>0$ is a small constant, and we may assume $r>0$, since
both $h_0$ and $h_1$ are even. From $s\G(s) = \G(s+1)$ and Stirling's
formula it follows that
$$
{\G(k - {5\over2}+iv)\over \G(k + {5\over2}+iv)} \ll k^{-5}
\qquad(12 \le k \le k_0).
$$
In the integral for $g(k)$ we shift the line of integration to $\R w = -2$,
taking $\R w = 2+\e$ as the line of integration in (2.6). Using Stirling's
formula and the above bound we obtain
$$
g(k) \ll QT^{-7/2}k^{-5},\leqno(6.2)
$$
and this bound can be further sharpened. Moreover directly from (6.1) we
have
$$
h_0(r) \;\ll\;Q{\roman e}^{-\pi r}.
\leqno(6.3)
$$
From (6.2) and (6.3) it is easily seen that the expressions in (4.5)
containing the functions $g(k)$ and $h_0(r)$ contribute $O(Q\log^6T)$.
It remains to deal with the contribution of $h_1(r)$. Since $h^*(w)$
is entire (to be rigorous, one has either to work with $h$ defined
by (1.6), or replace the constant 626 in (2.3) by a larger constant),
it transpires from (6.1) that in the expression for   $h_1(r)$
the poles of the integrand are at $w = \hf -n \;(n = 3,4,5,\ldots)$ and
at $w = m \pm ir\;(m = 0,-1,-2,\ldots)$.
The former ones are harmless and could be avoided by inserting factors
$r^2 + n^2 + {1\over4}$ in the numerator and denominator of $q(z)$ in (2.3).
We shift the line of integration
in the expression for   $h_1(r)$ to $\R w = -N$, letting eventually
$N \to\infty$. The main contribution will then come from the poles
at $w =\pm ir$ (these contributions are evaluated analogously), since
the residues at other poles are evaluated similarly, only they will be
of a lower order of magnitude. The residue at $w = -ir$ will be
$$
\ll |\G(2ir)|e^{2\pi r}|\G(\hf+ir)|^4|h^*(-ir)|
\ll e^{-\pi r}r^{-1/2}|h^*(-ir)|\leqno(6.4)
$$
with
$$
h^*(-ir) = \int_{\I z = -\e}zh(z){\G(-ir+iz)\over\G(1+ir+iz)}\d z,
$$
where $h^*$ is given by (2.6). Since $q(z) = 1 + O(|z|^{-4})$, it
is seen that $h^*(-ir)$ is majorized by two similar expressions, one
of which is ($z = T + Qy - i\e$)
$$
Q\left|\int_{-\infty}^\infty {T+Qy\over T+r+Qy}\,
{\roman e}^{-y^2+2i\e Q^{-1}y}\,{\G(iT-ir+iQy+\e)\over\G(iT+ir+iQy+\e)}\d y
\right|,\leqno(6.5)
$$
where we used $s\G(s) = \G(s+1)$. For $|y| \ge \log(rT)$ the portion of the
above integral is negligible, as is also the portion for
$r \ge T + T^\e Qy$, by Stirling's formula. Also note that
$|T+r+Qy| - |T-r+Qy| \le 2r$, so that the exponential function coming
from ${\roman e}^{-\pi r}$ in (6.4) and the gamma factors will have a
non-positive exponent. If
$$
T - T^\e Q \le r \le T + T^\e Q\leqno(6.6)
$$
holds, then from (2.2) and (4.8) we have
$$
h_1(r) \ll r^{-1/2}{\roman e}^{-\pi r}
\left|\int_{\Cal L}\G(-ir-iz)\G(-ir+iz)zh(z)\sinh(\pi z)\d z\right|,
\leqno(6.7)
$$
where $\Cal L$ is the real line with small indentations above and below
the points $z = -r$ and $z = r$, respectively.
It follows (by Stirling's formula) that the right-hand side of (6.7)
is of exponential decay if (6.6) holds.
Hence we are left with the most interesting range, namely
$$
1 \ll r \le T - T^\e Q.\leqno(6.8)
$$
Recall that the gamma-function admits an asymptotic expansion, for
$t \ge t_0 > 0$, whose first two terms
are
$$\eqalign{
\G(\s+it) &= \sqrt{2\pi}t^{-\s-{1\over2}}\exp\{-\hf\pi t +
i(t\log t - t + \hf\pi(\s-\hf))\}\cdot\cr&
\cdot\left(1 + \hf it^{-1}(\s - \s^2 -{\txt{1\over6}}) +O_\s(t^{-2})\right).
\cr}
$$
The quotient of gamma factors in (6.5) thus equals
$$
\left(1 + O\bigr({1\over T}\bigl)\right)
{\left({T-r+Qy\over T + r + Qy}\right)}^{\e-\hf}{\roman e}^{\pi r}
\exp(i\f(T,r,Q,y)),\leqno(6.9)
$$
where the term $O(1/T)$ admits an asymptotic expansion, and
by Taylor's formula we obtain
$$
\eqalign{
\f(T,r,Q,y) &= 2r + (T-r)\log(T-r) - (T+r)\log(T+r)\cr&
- 2Qy\left({r\over T} + {1\over3}\left({r\over T}\right)^3
+{1\over5}\left({r\over T}\right)^5 + \cdots\right)\cr&
+ {(Qy)^2\over T-r} + (T-r+Qy)\left(-{1\over2}{(Qy)^2\over (T-r)^2}
+ {1\over3}{(Qy)^3\over (T-r)^3} + \cdots\right)\cr&
- {(Qy)^2\over T+r} + (T+r+Qy)\left(-{1\over2}{(Qy)^2\over (T+r)^2}
+ {1\over3}{(Qy)^3\over (T+r)^3} + \cdots\right).\cr}\leqno(6.10)
$$
By (6.8) we have $Q|y|/(T\pm r) \le T^{-{1\over2}\e}$
for $|y| \le \log T$, so that we may truncate
the contribution of the last two  series above in such a way that
the tails will make a negligible contribution. The remaining terms
are inserted in
$$
\int_{-\log T}^{\log T}
{\roman e}^{-y^2+2i\e Q^{-1}y}\,{\G(iT-r+iQy+\e)\over\G(iT+r+iQy+\e)}\d y
= \int_{-\infty}^\infty + \;O\left({\roman e}^{-{1\over2}\log^2T}\right),
$$
where the term in (6.9) with the exponent $\e -\hf$ is again
simplified by Taylor's formula.
The integrals with the
remaining terms are evaluated by using the formula
$$
\int_{-\infty}^\infty y^j{\roman e}^{Ay-y^2}\d y
\= P_j(A){\roman e}^{{1\over4}A^2}
\qquad(j = 0,1,2,\ldots,\;P_0(A) = \sqrt{\pi}\,),
$$
where $P_j(z)$ is a polynomial in $z$ of degree $j$, which may be
explicitly evaluated by successive differentiation of the classic formula
$$
\int_{-\infty}^\infty {\roman e}^{Ay-y^2}\d y
\= \sqrt{\pi}{\roman e}^{{1\over4}A^2},
$$
considered as a function of $A$. The major contribution will come
from the term
$$
- 2Qy\left({r\over T} + {1\over3}\left({r\over T}\right)^3
+{1\over5}\left({r\over T}\right)^5 + \cdots\right)
$$
in $\f(T,r,Q,y)$, hence the total contribution will be, in view of (6.9),
$$
\ll {\roman e}^{-\pi r}r^{-1/2}|h^*(-ir)|
\ll Qr^{-1/2}\exp\left(-{Cr^2Q^2\over T^2}\right)\qquad(C>0).
$$
The analogous bound follows for the residue at $w = ir$.  In fact,
it follows that by the above procedure we obtain not only an upper bound,
but an asymptotic expansion
of the $h_1(r)$ in the range (6.8). This proves then the key bound (4.11),
establishes (5.7), and  completes the proof of Theorem 2.
\bigskip
\head
8. The asymptotic formula for sums of $H_j^3(\hf)$
\endhead
We shall present now the proof of the asymptotic formula (1.10) of Theorem
1. We start from (3.4)--(3.6), restricting ourselves as  to the range
$$
K^{\e} \;\le G \;\le K^{{1\over2}-\e},\leqno(7.1)
$$
and follow the approach developed in [5]. It is seen that it is the
term $\nu = 0$ in (3.4) whose contributions should be considered, because
the bound for the $\nu$-th term will be essentially the same as the
bound for the term $\nu = 0$, only it will be multiplied by $(G/K)^\nu$.
We note that the factors $\exp(-(f/K)^\lambda)$ and $U_\nu(fK)$ in
(2.1) can be conveniently removed by partial summation.
Next we follow the analysis carried out in [15, pp. 120 and 128-129] to
show that the contribution of $\nu = 3,5,6$ in (3.4) to (2.1) will be
small. Indeed, we have
$${\Cal H}_3(f;h_0)
\ll e^{-C\log^2K}\quad(C > 0)$$
and
$$
{\Cal H}_5(f;h_0) \ll d(f)f^{-1/2},\quad {\Cal H}_6(f;h_0)
\ll \s_{-1}(f)f^{1/2}K.
$$
The contribution of ${\Cal H}_4(f;h_\nu)$ was shown in
[5] to be $\ll GK^{1+\e}$. To estimate the contribution of
${\Cal H}_7(f;h_0)$ we note (see [3, Chapter 1]) that
$$
\sum_{n=1}^\infty \s_{2ir}(n)n^{-ir-s} \= \z(s-ir)\z(s+ir) \qquad
(r \in \RR,\, \R s > 1).
$$
Consequently by the Perron inversion formula (see e.g., [3, p. 486])
$$
\sum_{f\le3K}\,\s_{2ir}(f)f^{-{1\over2}-ir} \;\ll_\e\;
K^{2\mu({1\over2})+\e}
\qquad (K \ll |r| \ll K),\leqno(7.2)
$$
where as usual the Lindel\"of function $\mu(\s)$ is given by
$$
\mu(\s) = \limsup_{t\to\infty}
{\log|\z(\s+it)|\over\log t}.
$$
Instead of using directly (7.2) it is more expedient to use
the main contribution to the left-hand side of (7.2), which is
$$
{1\over2\pi i}\int_{\e-iU}^{\e+iU}\z(s+\hf-ir)\z(s+\hf+ir)K^s{\d s\over s}
\quad(K^\e \ll U \ll K^{1-\e}),
$$
and obtain a contribution which is, by the residue theorem,
$$\eqalign{&
\int_{\e-iU}^{\e+iU}
\int_{-\infty}^\infty{|\z(\hf+ir)|^4\over|\z(1+2ir)|^2}
h(r)\z(s+\hf-ir)\z(s+\hf+ir)K^s{\d s\over s}\d r\cr&
= \int_{-\infty}^\infty{|\z(\hf+ir)|^6\over|\z(1+2ir)|^2}h(r)\d r
+ \int_{-\e-iU}^{-\e+iU}
\int_{-\infty}^\infty\cdots \d s\d r + R\cr&
= J_1 + J_2 + R,\cr}
\leqno(7.3)
$$
say, where $R$ is the (small) contribution
from the integral over $[-\e\pm iU, \e\pm iU]$.
Alternatively, we may use the identity
$$
{\roman e}^{-Y^h} = {1\over2\pi i}\int_{(c)} Y^{-w}\G\bigl(1
+ {w\over h}\bigr){\d w\over w}\qquad(Y,h,c > 0)
$$
in (3.4) with $Y = f/K, h = C\log K$.

 After evaluating (3.1), we shall integrate it over $K$ from $K_0$
to $2K_0$, similarly as was done in Section 5. The integral $J_1$ in
(7.3) is the analogue of the integral on the left-hand side of (4.5).
Its total contribution will be $O(GK_0^{13/4}\log^{37/4}K_0)$, since
(4.14) holds and we use the best known estimate
$$
\int_0^T|\zt|^6\d t \ll T^{5/4}\log^{37/4} T,\leqno(7.4)
$$
which follows by H\"older's inequality from (4.2). The contribution
coming from $J_2$ will be analogous.
Namely note that the relevant range of $r$ in ${\Cal H}_7(f;h_0)$ is
$|r \pm K| \le G\log K$, hence it follows from (7.3) and the
argument given below that the total contribution
of ${\Cal H}_7(f;h_0)$ to the integrated version of (3.1) is
$$
\ll K_0^{3/2+\e}GU^{-1}(G + K_0^{2/3}) + GK_0^{13/4}\log^{37/4}K_0
$$
plus a quantity which is
$$\eqalign{
\ll \; &\int_{-U}^{U}\Bigl\{\int_{K_0}^{2K_0}
\int_{-\infty}^\infty K^{2-\e}\exp\Bigl(-(r-K)^2G^{-2}\Bigr)\log^2K_0
\times\cr&|\z(\hf+ir)|^4
|\z(\hf-\e+iu-ir)\z(\hf-\e+iu+ir)|\d r\d K\Bigr\}{\d u\over 1+|u|}.\cr}
\leqno(7.5)
$$
We shall take the maximum over $u$ in the integral in (7.5) and then
integrate; this will account for a loss of a log-factor in the final bound.
The integral in curly brackets resembles the one in (5.6), only it
has six and not eight zeta values, since now we are dealing with
$H_j^3(\hf)$ and not with $\H$. It equals $O(\exp(-c\log^2K_0))$ plus
$$
\eqalign{&
\int_{K_0-G\log K_0}^{2K_0+G\log K_0}|\z(\hf+ir)|^4|\z(\hf-\e+iu-ir)
\z(\hf-\e+iu+ir)|\times\cr&
\int_{K_0}^{2K_0}\exp\left(-(r-K)^2G^{-2}\right)
\d K\d r\cr&
= \left\{\int_{K_0+G\log K_0}^{2K_0-G\log K_0} +
\int_{K_0-G\log K_0}^{K_0+G\log K_0}
+ \int_{2K_0-G\log K_0}^{2K_0+G\log K_0}\right\}\cdots\d r\cr&
= I_1 + I_2 + I_3,\cr}
$$
say. The integrals $I_2$ and $I_3$ are estimated similarly. By H\"older's
inequality for integrals we have
$$
I_2 \ll G\Bigl(\int\limits_a^b|\z(\hf+ir)|^6\d r\Bigr)^{2\over3}
\Bigl(\int\limits_a^b|\z(\hf-\e+iu + ir)|^6\d r\Bigr)^{1\over6}
 \Bigl(\int\limits_a^b|\z(\hf-\e+iu - ir)|^6\d r\Bigr)^{1\over6}
$$
with $a = K_0 - G\log K_0$, $b = K_0 + G\log K_0$. Therefore we have
to estimate the integral of $|\zt|^6$ over a short interval. By using
the trivial estimate for $|\zt|^2$ and the asymptotic formula for the
integral of $|\zt|^4$ ([4, Chapter 5]) it follows that
$$
I_2 + I_3 \ll GK_0^{2\mu(\hf)+\e}(G + K_0^{2/3}).\leqno(7.6)
$$
There remains (on this occasion we fix $\e$)
$$
\eqalign{& I_1 =
G\int_{K_0+G\log K_0}^{2K_0-G\log K_0}|\z(\hf+ir)|^4|\z(\hf-\e+iu-ir)
\z(\hf-\e+iu+ir)|\times\cr&
\times\int_{(K_0-r)/G}^{(2K_0-r)/G}e^{-x^2}\d x\cdot\d r\cr&
= \sqrt{\pi}G\int_{K_0+G\log K_0}^{2K_0-G\log K_0}|\z|^4|\z||\z|\d r
+ O(\exp(-c\log^2K_0))\cr&
\ll
G\left(\int_{K_0}^{2K_0}|\z(\hf+ir)|^6\d r\right)^{2/3}
\left(\int_{K_0}^{2K_0}|\z(\hf-\e+iu+ ir)|^6\d r\right)^{1/6}\times\cr&
\times\left(\int_{K_0}^{2K_0}|\z(\hf-\e+iu - ir)|^6\d r\right)^{1/6}
\cr&
\ll GK_0^{{5\over4}+{\e\over3}}\log^{37/4}K_0,\cr}
$$
on using the functional equation for $\z(s)$ for the factors with "$-\e$" and
the bound (7.4). The gain of ${\e\over3}$ and one log-factor is more
than compensated by $K_0^{2-\e}\log^2K_0$ in (7.5).
We choose now  $U = K_0^{1/2-\e}$ and note that $\mu(\hf) < 1/6$ and
$G \le K_0^{1/2-\e}$.
It follows from (7.6) and the last bound that the total contribution of
${\Cal H}_7(f;h_0)$ to the integrated version of (3.1) is
$$
\ll GK_0^{13/4}\log^{37/4}K_0.\leqno(7.7)
$$
It remains to deal yet with the contribution of ${\Cal H}_2(f;h_0)$
and ${\Cal H}_1(f;h_0)$, which will produce the main term.
We have that the latter contributes
$$
4\pi^{-3/2}K^3G\left\{{\Cal C}^*_1(K,G) +
{\Cal C}^*_2(K,G)\right\} + O(K^{1+\e}G^3),\leqno(7.8)
$$
where
$$\eqalign{
{\Cal C}^*_1(K,G) &= \sum_{f\ge1}f^{-1}d(f)(\log K + \gamma -
\log(2\pi\sqrt{f})\exp(-(f/K)^\lambda),\cr
{\Cal C}^*_2(K,G) &=
-\sum_{f\ge1}f^{-1}d(f)(\log K + \gamma -
\log(2\pi\sqrt{f})U_0(fK),
\cr}
$$
and the function $U_0$ is given by (3.7). As in [15] we note that
${\Cal C}^*_1(K,G)$ equals
$${1\over2\pi i\lambda}\int_{(1)}
\left((\log K + \gamma - \log(2\pi))\z^2(w+1)
+ \z'(w+1)\z(w+1)\right)K^w\G(w/\lambda)\d w,
$$
and likewise ${\Cal C}^*_2(K,G)$ can be represented by a similar type of
integral. The line of integration is shifted to $\R w = -1$,
where the integrand is regular. There is a pole of order three at
$w = 0$, hence by the residue theorem and
Stirling's formula for $\G(s)$ we obtain
$$
\eqalign{
{\Cal C}^*_1(K,G) &= \sum_{j=0}^3A_j\log^jK + O(K^{\e-1}),\cr
{\Cal C}^*_2(K,G) &= \sum_{j=0}^3B_j\log^jK + O(K^{\e-1}),\cr}\leqno(7.9)
$$
with $A_3 = B_3 = 1/3$. The $O$-term in (7.8) comes from the fact
(see the definition of ${\Cal H}_1(f;h)$ in (3.6))  that we have
$$
\eqalign{
({\hat h}_0)'(\hf) &= 2i\pi^{3/2}K^3G + O(KG^3),\cr
({\hat h}_0)''(\hf) &= 8i\pi^{3/2}K^3G\log K + O(KG^3\log K).\cr}
$$

\medskip
From (7.3)--(7.9) we obtain ($G = G(K_0) \,(\le K_0^{1/2-\e})$
will be suitably chosen a  little later; see (8.9))
$$\eqalign{
\int_{K_0}^{2K_0}{\Cal C}(K,G)\d K  &=
GK^4{\bar P}_3(\log K)\Bigg|_{K_0}^{2K_0} \cr&
+ O(GK_0^{13/4}\log^{37/4}K_0)
+ O(G^3K_0^{2+\e}),\cr}
\leqno(7.10)
$$
where ${\bar P}_3$ is another cubic polynomial, this time with leading
coefficient  $2/(3\pi^{3/2})$. Here we have assumed that the
total contribution of ${\Cal H}_2(f;h)$ can be absorbed in the error
terms in (7.10), which will be shown in Section 8 with suitable $G$.

On the other hand, applying (1.7) in the form
$$
\sum_{K\le\k_j\le K+H}\a_j H_j^3(\hf) \ll_\e K^{1+\e}H
\qquad(1 \ll H \le K)
$$
and using the method of proof of Section 5, it is seen that
$$\eqalign{&
\int_{K_0}^{2K_0}{\Cal C}(K,G)\d K
= \sum_{j\ge1}\a_jH_j^3(\hf)\int_{K_0}^{2K_0}(\k_j^2 + {\txt{1\over4}})
\exp(-(\k_j-K)^2G^{-2})\d K + o(1)\cr&
= \sqrt{\pi}G\sum_{K_0\le\k_j\le2K_0}\a_jH_j^3(\hf)\k_j^2  +
O(K_0^{3+\e}G^2).
\cr}\leqno(7.11)
$$
Therefore we obtain from (7.10) and (7.11)
$$\eqalign{&
\sum_{K_0\le\k_j\le2K_0}\a_jH_j^3(\hf)\k_j^2
= K^4\left({2\over3\pi^2}\log^3K + a_2\log^2K + a_1\log K + a_0\right)
\Bigg|_{K_0}^{2K_0}\cr&
+ O(K_0^{13/4}\log^{37/4}K_0) + O(GK_0^{3+\e})  \cr}\leqno(7.12)
$$
plus the contribution of ${\Cal H}_2(f;h)$. We apply partial summation
(to get rid of $\k_j^2$), replace $K_0$ by $K_02^{-j}$, and sum over $j$.
The $O$-terms will be absorbed in the  $O$-term of Theorem 1
if $G = K_0^\a$ with any $0 < \a < 1/4$.

\head
8. The contribution of ${\Cal H}_2(f;h)$
\endhead

To complete the proof of Theorem 1 it remains to show that the total
contribution of ${\Cal H}_2(f;h)$ is absorbed in the $O$-terms in (7.12)
with suitable $G$. We follow, as before, the proof given in [5].
We use the observation made in [7] which states that the relevant sum to be
estimated is, after integration over $\,[K_0,\,2K_0]\,$,
$$
\eqalign{&
GK_0^{5/2}\sum_{f\le 3K_0}f^{-1/2}\sum_{m\le fG^{-2}\log^2K_0}(m/f)^{1/4}
d(m)d(m+f)\times\cr&
\times
\left(\sqrt{m\over f} + \sqrt{1+ {m\over f}}\right)^{-2iK_0}
{\roman e}^{-CG^2mf^{-1}}\log
\left(\sqrt{{m\over f}} + \sqrt{1+ {m\over f}}\,\right)^{-1}.\cr}\leqno(8.1)
$$
Note that (8.1) corresponds to (3.1) of [5] with the additional factor
$(m/f)^{1/4}$, namely to (16) of [7]. As in (3.2) of [5] we replace
$m+f$ by $n$ and consider subsums of the sum in (8.1) where $m \sim M$
(meaning $M < m \le 2M$), $n \sim N$. If we get rid of the last two
factors in (8.1) by partial summation and Taylor's formula, respectively,
we are left with the sum
$$
\eqalign{&
GK_0^{5/2}\sum_{n\sim N}d(n)n^{-1/4}\sum_{m\sim M}d(m)m^{-3/4}
\exp(iF(m,n)),\cr&
F(m,n) \;:=\;-2K_0\log\left(\sqrt{m\over n-m} + \sqrt{n\over n-m}\,\right),
\cr}\leqno(8.2)
$$
and we have, with effectively computable constants $b_j$,
$$
\log\left(\sqrt{m\over n-m} + \sqrt{n\over n-m}\,\right)
= \sum_{j=1}^\infty b_j{\left({m\over n}\right)}^{j/2}.\leqno(8.3)
$$
As in [5, eq. (3.4)], we have the conditions
$$
K_0^\e \le G \le K_0^{1/2-\e},\; MG^2\log^2K_0 \ll N \ll K_0.\leqno(8.4)
$$
By applying the Cauchy-Schwarz inequality we see that the sum in (8.2) is
$$
\eqalign{&
\le\left(\sum_{n\sim N}d^2(n)n^{-1/2}\right)^{1/2}
\left(\sum_{n\sim N}\left|\sum_{m\sim M}d(m)m^{-3/4}\exp(iF(m,n))\right|^2
\right)^{1/2}\cr&
\ll N^{1/4}\log^2N{\sum\limits}^{1/2},\cr}
$$
where we have set
$$
\eqalign{& \sum
\;:=\; \sum_{n\sim N}\left|\sum_{m\sim M}d(m)m^{-3/4}\exp(iF(m,n))\right|^2
\cr&
= \sum_{m\sim M}d^2(m)m^{-3/2}O(N) \cr&
+ \sum_{m_1\not=m_2}d(m_1)d(m_2)(m_1m_2)^{-3/4}\sum_{n\sim N}
\exp(iF(m_1,n)-iF(m_2,n))\cr&
\ll NM^{\e-1/2} + M^{\e-3/2}\sum_{m_1\not=m_2}
\left|\sum_{n\sim N}\exp(iF(m_1,n)-iF(m_2,n))\right|.\cr}
$$
The effect of this procedure is that the exponential sum over $n$
does not contain the divisor function, and consequently can be
estimated by the technique of exponent pairs (see e.g., [3, Chapter 2]).
Note that by (8.3) we have (in the relevant range for $m,n$)
$$
{\partial \over \partial n}(F(m_1,n)-F(m_2,n)) \asymp
|m_1 - m_2|K_0M^{-1/2}N^{-3/2}.
$$
Thus if $(\k,\lambda)$ is an exponent pair, then we have
$$
\eqalign{&
\sum \ll NM^{\e-1/2} + M^{\e-3/2}\sum_{m_1\not=m_2}
\left({N^{3/2}M^{1/2}\over |m_1-m_2|K_0} +
\left({K_0M^{1/2}\over N^{3/2}}\right)^\k
N^\lambda\right)\cr&
\ll NM^{\e-1/2} + N^{3/2}K_0^{\e-1} + M^{{1\over2}+{\k\over2}}K_0^\k
N^{\lambda-{3\over2}\k}.
\cr}
$$
Hence in view of (8.4) the expression in (8.2) is bounded by
$$
\eqalign{&
GK_0^{{5\over2}+\e}\left(N^{3\over4}M^{-{1\over4}} + NK_0^{-{1\over2}}
+ K_0^{\k\over2}M^{{1\over4}+{\k\over4}}N^{{\lambda\over2}+
{1\over4}-{3\over4}\k}\right)\cr&
\ll GK_0^{{5\over2}+\e}N^{3\over4}M^{-{1\over4}}+ GK_0^{3+\e} +
GK_0^{{11\over4}+\e}M^{3\over8}N^{1\over8}\cr&
\ll GK_0^{{5\over2}+\e}N^{3\over4}M^{-{1\over4}} + GK_0^{3+\e}
+ GK_0^{{5\over2}+\e}K_0^{1\over4}(NG^{-2})^{3\over8}N^{1\over8}\cr&
\ll GK_0^{{5\over2}+\e}N^{3\over4}M^{-{1\over4}} + GK_0^{3+\e}
+ GK_0^{{13\over4}+\e}G^{-{3\over4}}
\cr}     \leqno(8.5)
$$
with $(\k,\,\lambda) = (\hf,\,\hf)$.
The bound in (8.5) will be used for large $M$. For small $M$ we shall
transform the sum
$$
S(N) := \sum_{{1\over2}N\le n\le{5\over2}N}\f(n)d(n)n^{-1/4}\exp(iF(m,n))
$$
by Voronoi's summation formula (see e.g., [3, Chapter 3]), treating the
real and imaginary part separately. Here $\f(x) \ge 0$ is a smooth
function supported in $[\hf N,\,{5\over2}N]$ such that it equals unity
in $[N,\,2N]$ and $\f^{(r)}(x) \ll_r N^{-r}\,(r = 0,1,\ldots\,)$.
Then we have
$$
\eqalign{
S(N) &= \int_{{1\over2}N}^{{5\over2}N}(\log x + 2\gamma)x^{-1/4}\f(x)
\exp(iF(m,x))\d x \cr&
+ \sum_{n=1}^\infty \int_{{1\over2}N}^{{5\over2}N}
\f(x)x^{-1/4}\a(nx)  \exp(iF(m,x))\d x,\cr}\leqno(8.6)
$$
where $\a(nx)$ admits an asymptotic expansion whose first term is
$$
-2^{1/2}(xn)^{-1/4}\sin(4\pi\sqrt{nx} - \pi/4).
$$
By the first derivative test the first integral in (8.6) is
$$
\ll {N^{5/4}\log N\over M^{1/2}K_0},
$$
hence it contributes to (8.2)
$$
\ll GK_0^{3/2}N^{5/4}M^{-1/4}\log^2K_0 \ll GK_0^{11/4}\log^2K_0.
$$
Further consider the main contribution of the terms in (8.6),
which is a multiple of
$$
\int_{{1\over2}N}^{{5\over2}N}\f(x)x^{-1/2}n^{-1/4}
\exp\left(4\pi i\sqrt{nx}
\pm iK_0\sum_{j=1}^\infty b_j\left({m\over x}\right)^{\hf j} \right)\d x.
\leqno(8.7)
$$
The case of the ``minus" sign is less difficult, and in the case of the
``plus" sign, let
$$
f(x) =  f(x;m,n,K_0) :=  4\pi \sqrt{nx} +
K_0\sum_{j=1}^\infty b_j\left({m\over x}\right)^{\hf j},
$$
so that
$$
{\partial f\over\partial x} = 2\pi \sqrt{n\over x}
- K_0 \sum_{j=1}^\infty \hf jb_jm^{j/2}x^{-j/2-1}.
$$
If $n > CK_0^2MN^{-2}$ with sufficiently large $C > 0$,
then $ {\partial f\over\partial x}
\asymp \sqrt{n\over x}$. Therefore the above integral becomes,
on integrating by parts,
$$
in^{-1/4}\int_{{1\over2}N}^{{5\over2}N}
{\left({\f(x)x^{-1/2}\over {\partial f\over\partial x}}\right)}'
\exp(if(x))\d x.
$$
But as
$$
{\left({\f(x)x^{-1/2}\over {\partial f\over\partial x}}\right)}'
\ll {1\over\sqrt{nxN}},
$$
it follows by repeated integration by parts that the contribution
of  $n > CK_0^2MN^{-2}$ is negligible. If $n \le CK_0^2MN^{-2}$, then
the exponential integral in question may have a saddle point $x_0$,
namely the solution of ${\partial f\over\partial x}=0$. Hence
$$
2\pi\sqrt{n\over x_0} = K_0\sum_{j=1}^\infty \hf jb_jm^{j/2}x_0^{-j/2-1},
$$
giving (since $b_1 = 1$)
$$
x_0 \,\sim\,{K_0\over2\pi}\sqrt{m\over n},
$$
and $x_0 \in [\hf N,\,{5\over2}N]$ for $n \asymp K_0^2MN^{-2}$.
By the saddle point method (see [3, Chapter 2]) the main contribution
comes from the saddle point and is
$$
\ll \left|{\partial^2 f\over\partial x^2}\Big|_{x=x_0}\right|^{-1/2}
\ll \left({1\over N}\sqrt{n\over N}\right)^{-1/2} = N^{3/4}n^{-1/4}.
$$
Thus the integral in (8.7) is
$\ll N^{1/4}n^{-1/2}$, and consequently the sum in (8.6) is
$$
\ll N^{1/4}\sum_{n\le CK_0^2MN^{-2}}d(n)n^{-1/2}
\ll K_0M^{1/2}N^{-3/4}\log K_0,
$$
and the total contribution is therefore
$$
\ll GK_0^{11/4}\log^2K_0 + GK_0^{7/2}M^{3/4}N^{-3/4}\log^2 K_0. \leqno(8.8)
$$
Hence for $M \ge N^{3/2}/K_0$ we use (8.5) and otherwise we apply (8.8); if
$N \le K_0^{2/3}$ then $N^{3/2}/K_0 \le 1$, but then we can simply use
(8.5). We obtain, in view of (7.11) and (7.12) and the discussion
thereafter, that the total contribution
of the error terms in Theorem 1 will be
$$
\ll K_0^{5/4}\log^{37/4}K_0 + GK_0^{1+\e} +
K_0^{{5/4}+\e}G^{-3/4} \ll  K_0^{5/4}\log^{37/4}K_0
$$
for
$$
G = K_0^{1/7}.\leqno(8.9)
$$
This completes the proof of Theorem 1.
Note that, apart from the contribution of the integral with six zeta
values (cf. (7.3)), the remaining terms are of the order $K_0^{8/7+\e}$
with the choice $G = K_0^{1/7}$, and more refined exponential sum
techniques could yield even smaller values of $G$.
From (7.12) it follows that the leading coefficient of $P_3(x)$ in
(1.10) is $4/(3\pi^2)$.

\bigskip

\head
9. Another proof of Theorem 2
\endhead

We shall sketch now another proof of Theorem 2 (cf. (4.1)), namely
$$
\sum_{\k _j \leq K}\alpha _jH_j^4(\hf ) + O\left (\log ^2K\int_0^K
|\zt|^8\d t \right ) =K^2P_6(\log K) +
O(K^{4/3+\e }).\tag9.1
$$
The argument is based on M. Jutila's proof [7] of (1.9), and will be
outlined below. Similarly as in the proof of Theorem 1, it
is the contribution of ${\Cal H}_2(f;h)$ (see (3.6)) that is
the essential one. To introduce $H^2_j(\hf)$ in Motohashi's
transformation formula for sums of $H^2_j(\hf)$ ([15, Lemma 3.8])
and obtain the formula for sums of $\H$, one uses [7, Lemma 1].
This formula says that
$$\eqalign{
H_j^2(\hf) &= \sum_{mn\le 3K^2}t_j(m)t_j(n)(mn)^{-1/2}
\exp(-(mn/K^2)^\lambda)\cr&
-   \sum_{mn\le 3K^2}t_j(m)t_j(n)(mn)^{-1/2}R_j(mnK^2) + O(1),\cr}
\tag9.2
$$
for $|\k_j - K| \le G\log K$ with $\log^2K < G < K^{1-\delta}$ for
$0 < \delta < 1$, $\lambda = C\log K$ with sufficiently large $C>0$.
The function $R_j$ in (9.2) comes from the squaring of the functional
equation for $H_j(\hf + w)$, namely
$$
\eqalign{
R_j(x) &= {1\over2\pi^4i\lambda}\int_{-\lambda^{-1}-i\lambda^2}
^{-\lambda^{-1}+i\lambda^2}(16\pi^4x)^w\G^2(\hf-w+i\k_j)\G^2(\hf-w-i\k_j)
\times\cr&
\times(\cosh(\pi\k_j) + \sin(\pi w))^2\G(w/\lambda)\d w.\cr}
$$
In the context of [7] the error term $O(1)$ in (9.2) suffices, but
similarly to [15, Lemma 3.9] this error term can be considerably sharpened.
The main term (i.e., $K^2P_6(\log K)$ in (9.1)) is derived analogously
as was done in the proof of Theorem 1; it is obtained
in terms of the expressions resembling the functions
${\Cal C}^*_j\;(j=1,2)$ in (7.8), only in this case they
will be somewhat more complicated. Namely to obtain the asymptotic
formula for the sum
$$
\sum_{j=1}^\infty \a_jH_j^4(\hf)h_0(\k_j) \leqno(9.3)
$$
with $h_0$ given by (1.6), we use the Mellin relation
$$
\exp(-x^\lambda) = {1\over2\pi i\lambda}\int_{(1)}\G(z/\lambda)x^{-z}\d z
\qquad(x,\lambda >     0)
$$
in conjunction with (9.2) and [15, Lemma 3.8]. We use the
identity (3.3) to transform the product of two $t_j$-functions
into one, and extend summation over all values of $m,n$, producing
a negligible error. Then we obtain two divisor functions, and we use
the classical identity
$$
\sum_{n=1}^\infty d^2(n)n^{-s} \;=\; {\z^4(s)\over\z(2s)}\quad(\R s > 1).
$$
It follows that, similarly to the case
of Theorem 1, the main term for (9.3) will be of the form
$$
4\pi^{-3/2}K^3G({\Cal D}_1^*(K,G) +   {\Cal D}_2^*(K,G)),
$$
where   ${\Cal D}_1^*(K,G)$ comes from the first sum on the right-hand
side of (9.2). We have   ($\gamma$ is Euler's constant)
$$\eqalign{
{\Cal D}_1^*(K,G) &=  {1\over2\pi i\lambda}\int_{(1)}
\Bigl\{(\log K + \gamma - \log(2\pi)){\zeta^4(w+1)\over\zeta(2w+2)}\cr&
+ {1\over2}{\left({\zeta^4(w+1)\over\zeta(2w+2)}\right)}'\Bigr\}
\z(2w+1)K^{2w}\G(w/\lambda)\d w,\cr}\leqno(9.4)
$$
and analogously   ${\Cal D}_2^*(K,G)$ comes from the   second sum
on the right-hand side of (9.2). The integrand in (9.4) has a pole
of order six at $w=0$. We shift the line of integration to $\R w = -1$,
developing the integrand into power series to calculate the residue.
The coefficient of $\log^6K$ is found to be $4/(15\pi^2)$, and clearly
the coefficients of lower powers of the logarithm can be also
evaluated explicitly. This is the analogue of $A_3 = 1/3$ in (7.9).
The  coefficient of $\log^6K$ coming from ${\Cal D}_2^*(K,G)$
will be the same.
Proceeding as was done in Section 7, we see then that the leading
coefficient of $P_6(x)$ in (1.11) is $16/(15\pi^4)$, as claimed.

\medskip

We continue now the second proof of Theorem 2.
From the discussion above it is seen that the relevant sum to be estimated
(this corresponds to [7, eq. (16)]) is, up to a constant factor,
$$\eqalign{&
GK^{5/2}\sum_{f\ll K^2} v(f)d(f)f^{-3/4}\sum_{m \leq fG^{-2}\log ^2 K}
m^{-1/4}d(m)d(m+f)\cr&
\times \left (\sqrt{\frac mf }+\sqrt{1 +\frac mf }\,\right )^{-2iK}
\exp \left (-G^2\log ^2
\left (\sqrt{\frac mf }+\sqrt{1 +\frac mf }\,\right )\right ) \cr&
\times \left (\log
\left (\sqrt{\frac mf }+\sqrt{1 +\frac mf }\,\right )\right )^{-1},\cr}
$$
where $v$ is a smooth weight function supported in $\,[F,\,2F]\,$ with
$F \ll K_0^2$, and $K_0 \le K \le 2K_0$. A new ingredient is the
last log-factor (coming from integration), which is of the order
$\ll \sqrt{f/m}$. Consider now the sum over $f\asymp F$ and
$m \asymp M$. Then, by the above remarks,  the final estimate
in \cite{7}, namely
$$
\ll GK^{\e}(F^{-1/2}KM^{1/2})^{3/2},
$$
should be modified by cancelling the factor $G$ and
multiplying by $\sqrt{F/M}$. Therefore the contribution coming from
${\Cal H}_2(f;h)$ will be
$$
\ll K_0^{3/2+\e }(M/F)^{1/4} \ll K_0^{4/3+\e},
$$
since $M/F \ll G^{-2}\log ^2 K_0$ and $G = K_0^{1/3} \,(\asymp Q\,$ of
Section 4). This finishes the discussion concerning
the second proof of Theorem 2.

\bigskip
\head
10. The first moment of $H_j(\hf)$
\endhead

As promised in the Introduction, we shall say a few words at the end
on the sum
$$
\sum_{\k_j\le T}\a_j H_j(\hf).\leqno(10.1)
$$
In conjunction with the conjecture (1.12) I expect the sum in (10.1)
to be equal to
$$
AT^2 + O(T\log^3T)\qquad(A = {1\over\pi^2}),\leqno(10.2)
$$
where the error term in (10.2) comes from the integral
with $|\z(\hf + it)|^2$ in (1.12), and the value of $A$ is provided
by Random matrix theory (see the discussion at the end of Section 1).
However obtaining (10.2) is rather
difficult. Namely, simple specialization
(simplification) of the procedure used by Y. Motohashi [14]
for sums of $H^2_j(\hf)$ does not work directly. In any case it can
be shown that
$$
T^2(\log T)^{-7/2} \ll \sum_{\k_j \le T}\a_jH_j(\hf) \ll T^2(\log T)^{1/2}.
\leqno(10.3)
$$
The upper bound in (10.3) follows from the Cauchy-Schwarz inequality
and (1.2). To derive the lower bound, let
$$
S(T) := \sum_{T\le \k_j \le 2T}\a_jH_j(\hf).
$$
For a given $V > 0$ we have (since $H_j(\hf) \ge 0$)
$$
S(T) \ge V \sum_{T\le \k_j \le 2T,H_j(\hf)\ge V}\a_j,
$$
and we obtain
$$\eqalign{
T^2\log T &\ll \sum_{T\le \k_j \le 2T}\a_jH_j^2(\hf)
= \sum_{H_j(\hf) \ge V} + \sum_{H_j(\hf) < V}\cr&
\ll \left( \sum_{T\le \k_j \le 2T,H_j(\hf)\ge V}\a_j
\sum_{T\le \k_j \le 2T}\a_jH_j^4(\hf)\right)^{1/2} +
V^2\sum_{T\le \k_j \le 2T}\a_j\cr&
\ll \left(V^{-1}S(T)T^2\log^6T\right)^{1/2} + T^2V^2.\cr}
$$
Here we used the best possible bounds (cf. [4, eq. (5.48)] and (1.11))
$$
\sum_{\k_j \le T}\a_j \ll T^2,\quad
\sum_{\k_j \le T}\a_j H_j^4(\hf) \ll T^2\log^6T.
$$
The choice $V = \delta\sqrt{\log T}$ for sufficiently small $\delta > 0$
yields then
$$
T^4\log^2T \ll V^{-1}S(T)T^2\log^6T,
$$
giving the lower bound in (10.3).

\smallskip
One way to tackle the sum in (10.1) is
to take $ n=1$ in Kuznetsov's trace formula
([14, eq. (2.5)]) and multiply by $m^{-u}$ to obtain
$$
\eqalign{&
\sum_{j=1}^\infty \e_j\a_j t_j(m)m^{-u}h(\k_j)\cr&
= -{1\over\pi}\int_{-\infty}^\infty \s_{2ir}(m)m^{-u-ir}
{h(r)\over|\z(1+2ir)|^2}\d r + \sum_{\ell=1}^\infty m^{-u}
\ell^{-1}S(m,-1;\ell)\psi\bigl(4\pi{\sqrt{m}\over\ell}\bigr),\cr}\leqno(10.4)
$$
where $S(m,n;\ell)$ is the Kloosterman sum, $h(r)$ is given
by (2.3), while   with $h^*(s)$ given by (2.6) we set
$$
\psi(x) = \frac{1}{\pi^2}\int_{(\a)}\frac{(x/2)^{-2s}}{\cos(\pi s)}
h^*(s)\d s\qquad(-3/2 < \a < 3/2).
\leqno(10.5)
$$

\medskip
We proceed now,  assuming that $\R u > 2$
and $\a = -2/3$ in (10.5). Using the trivial bound $|S(m,-1;\ell)| \le \ell$,
we note that summation over $m$ in (10.4) yields, by absolute convergence,
$$
\eqalign{&
\sum_{j=1}^\infty \e_j\a_j H_j(u)h(\k_j) +
{1\over\pi}\int_{-\infty}^\infty \z(u+ir)\z(u-ir)
{h(r)\over|\z(1+2ir)|^2}\d r\cr&
= \sum_{m=1}^\infty m^{-u}\sum_{\ell=1}^\infty\ell^{-1}S(m,-1;\ell)
\psi(4\pi{\sqrt{m}\over\ell}).\cr}\leqno(10.6)
$$
By deforming suitably the contour and applying the residue theorem, we
see that the integrated term admits analytic continuation to the
region $\R u < 1$ which is of the form
$$
{1\over\pi}\int_{-\infty}^\infty \z(u+ir)\z(u-ir)
{h(r)\over|\z(1+2ir)|^2}\d r + 4{h(i(u-1))\over\z(3-2u)}.
$$
Since $H_j(\hf) = 0$ if $\e_j = -1$ and $h(\pm\hf i) = 0$,  (10.6)
reduces to  (compare with (1.12) when $k=1$)
$$
\sum_{j=1}^\infty \a_jH_j(\hf)h(\k_j) +
{1\over\pi}\int_{-\infty}^\infty |\z(\hf+ir)|^2
{h(r)\over|\z(1+2ir)|^2}\d r = L(\hf),\leqno(10.7)
$$
where $L(u)$ is the analytic continuation of the function
$$
\sum_{m=1}^\infty m^{-u}\sum_{\ell=1}^\infty\ell^{-1}S(m,-1;\ell)
\psi(4\pi{\sqrt{m}\over\ell})\qquad(\R u > 2).\leqno(10.8)
$$
One can try to transform the expression for $L(u)$ by using
the properties of the Kloosterman--Selberg zeta-function
$$
Z_{m,n}(s) := (2\pi\sqrt{mn})^{2s-1}\sum_{\ell=1}^\infty
S(m,n;\ell)\ell^{-2s}\qquad(\R s > 1).
$$
Namely one has the spectral decomposition (see [4, eqs. (5.65)--(5.68)]
of  $Z_{m,n}(s)$. This can be used in (10.8), and one
expects that the main contribution will come
from the discrete spectrum (i.e. [4, (5.66)]). However this will lead
eventually to the same type of sum as the one we started from.

\smallskip
One can follow the approach of [14] and write ($-3/2 < \a < -1/4$)
$$
L(u) = \pi^{-2}\sum_{\ell=1}^\infty\ell^{-1}P(u;\ell),
$$
$$
P(u;\ell) = \int_{(\a)}(2\pi/\ell)^{-2s}{h^*(s)\over\cos(\pi s)}
Q(s;u,\ell)\d s,
$$
$$
\eqalign{&
Q(s;u,\ell) = \sum_{m=1}^\infty m^{-u-s}S(m,-1;\ell)\cr&
= \sum_{(a,\ell)=1,a{\bar a}\equiv 1(\roman{mod}\ell)}{\roman e}(-a/\ell)
E(u+s;{\roman e}({\bar a}/\ell)),\cr}\leqno(10.9)
$$
where $E$ is the Lerch zeta-function ($1 \le h \le k,\,k\ge 2,
\,h,k\in\NN,\;{\roman e}(z) = {\roman e}^{2\pi i z}$)
$$
E\left(s;{\roman e}\left({h\over k}\right)\right) := \sum_{m=1}^\infty
{\roman e}\left({mh\over k}\right)m^{-s}
= \sum_{j=1}^k{\roman e}\left({jh\over k}\right)
k^{-s}\z(s,{j\over k}),
$$
initially defined for $\R s > 1$. It
can be expressed in terms of the Hurwitz zeta-function,
defined for $0 < a \le 1,\,\s>1$ by
$\z(s,a) = \sum_{n=0}^\infty (n+a)^{-s}$. Since
$\z(s,{j\over k})$ has a only the simple pole at $s=1$ with residue 1,
it follows that $E$  is entire, and satisfies the functional equation
$$
E\left(s;{\roman e}\left({h\over k}\right)\right)
= {\G(1-s)\over(2\pi)^{1-s}}\left\{{\roman e}^{{\pi i\over2}(1-s)}
\z(1-s,{h\over k}) + {\roman e}^{{\pi i\over2}(s-1)}
\z(1-s,1 - {h\over k})\right\}.  \leqno(10.10)
$$
This means that the second expression in (10.9) provides the analytic
continuation of $Q(s;u,\ell)$ as an entire function of both $u$ and $s$,
of polynomial growth in $|u|+|s|$.

\smallskip
This, however, differs from Motohashi's situation [14],
where he obtained the Estermann zeta-function $D$, represented in the
region of absolute convergence by the series
$$
D(s,\xi;{\roman e}(b/\ell)) \;:=\;
\sum_{n=1}^\infty n^{-s}\s_\xi(n){\roman e}(nb/\ell)
\qquad(1 \le b \le \ell;\, b,\ell\in\NN).
$$
This function has two simple poles (at $s = 1$ and $1+\xi$)
which are (in part) responsible for the main term $(2.34)_1$   in [14].
But we do not have such a term here!  What we get is simply,
since $E$ is entire,
$$
L(\hf) = {1\over\pi^2}\int_{(\a)}(2\pi)^{-2s}{h^*(s)\over\cos(\pi s)}
\sum_{\ell=1}^\infty\sum_{(a,\ell)=1}{\roman e}(-a/\ell)\ell^{2s-1}
E(s+\hf;{\roman e}({\bar a}/\ell))\d s.\leqno(10.11)
$$
In (10.11) we have $-3/2 < \a < -1/2$.
To transform further $L(\hf)$ we make
the change of variable $s = \hf - w$ in (10.11) and use the functional
equation (10.10). It follows that $L(\hf)$ is a linear combination of
$$
I_+ := \int_{(\b)}(2\pi)^{w}h^*(\hf-w){\G(w)\over\sin(\hf \pi w)}
M_+(w)\d w\quad(1 < \b < 2)
$$
and
$$
I_- := \int_{(\b)}(2\pi)^{w}h^*(\hf-w){\G(w)\over\cos(\hf \pi w)}
M_-(w)\d w\quad(1 < \b < 2),
$$
where for $\R w > 1$
$$
\eqalign{
M_+(w) &:= \sum_{\ell=1}^\infty \sum_{(a,\ell)=1,a{\bar a}\equiv 1
(\roman{mod}\ell)}{\roman e}(-a/\ell)
l^{-2w}\left((\z(w,{\bar a}/\ell) + \z(w, 1 - {\bar a}/\ell)\right),\cr
M_-(w) &:= \sum_{\ell=1}^\infty \sum_{(a,\ell)=1,a{\bar a}
\equiv 1(\roman{mod}\ell)}{\roman e}(-a/\ell)
l^{-2w}\left((\z(w,{\bar a}/\ell) - \z(w, 1 - {\bar a}/\ell)\right).\cr}
$$
The problem is to obtain analytic continuation of the functions
$M_\pm(w)$ to the left of the line $\R w = 1$, since one would like to
move the contour of integration in $I_+$ and $I_-$ to the left.

\medskip

It transpires that in any case it seems difficult to show that
the sum in (10.1) equals the expression in (10.2).

\vfill
\eject\topskip2cm
\Refs
\bigskip

\item{[1]} J.B. Conrey, $L$-functions and random matrices, in ``Mathematics
Unlimited" (Part I), B. Engquist and W. Schmid eds., Springer,
2001, pp. 331-352.

\item{[2]} J.B. Conrey, D.W. Farmer, J.P. Keating, M.O. Rubinstein
and N.C. Snaith, Integral moments of $L$-functions, preprint, 58pp,
arXiv:math.NT/0206018,

{\tt http://front.math.ucdavis.edu/mat.NT/0206018}.

\item{[3]} A. Ivi\'c, The Riemann zeta-function, John Wiley \&
Sons, New York, 1985.

\item{[4]} A. Ivi\'c,  The mean values of the Riemann zeta-function, Tata
        Institute of Fundamental Research, Lecture Notes {\bf82},
    Bombay 1991 (distr. Springer Verlag, Berlin etc.), 363 pp.

\item{[5]} A. Ivi\'c,
On sums of Hecke series in short intervals, J. de Th\'eorie des
Nombres Bordeaux {\bf14}(2001), 554-568.

\item{[6]} H. Iwaniec, Small eigenvalues of the Laplacian for $\,\G_0(N)\,$,
Acta Arithmetica {\bf56}(1990), 65-62.

\item{[7]} M. Jutila, The fourth moment of central values of Hecke
series, in Number Theory, Proc. of the Turku Symposium 1999,
Walter de Gruyter, Berlin, 2001, 167-177.

\item{[8]} M. Jutila and Y. Motohashi, A note on the mean value
of the zeta and $L$-functions XI, Proc. Japan Acad. {\bf78}, Ser. A
(2002), 1-6.

\item {[9]} S. Katok and P. Sarnak, Heegner points, cycles and Maass
forms,  Israel J. Math. {\bf84}(1993), 193-227.

\item{[10]} N.V. Kuznetsov, Sums of Kloosterman sums and the eighth
power moment of the Riemann zeta-function, T.I.F.R. Stud. Math.
{\bf13}(1989), 57-117.

\item{[11]} N.V. Kuznetsov, The true order of the Riemann
zeta-function on the critical line (preprint), Institute for Applied
Math., Khabarovsk, 1998,  88pp.

\item{[12]} N.V. Kuznetsov, The Hecke series at the center of the
critical strip (preprint, in Russian), Vladivostok: Dal'nauka, 1999, 27pp.

\item{[13]} Y. Motohashi, Kuznetsov's paper on the eighth power
moment of the Riemann zeta-function (revised, Part I),
manuscript dated June 22, 1991, 24pp.

\item{[14]} Y. Motohashi, Spectral mean values of Maass wave form
$L$-functions, J. Number Theory {\bf42}(1992), 258-284.

\item{[15]} Y. Motohashi, Spectral  theory of the Riemann zeta-function,
Cambridge University Press, 1997.

\item{[16]} Y. Motohashi, A functional equation for the spectral fourth
moment of modular Hecke $L$-functions,
Proceedings of the conference ``Workshop in
Analytic Number Theory", June 24-28, 2002 at the Max Planck Institut
f\"ur Mathematik, Bonn, 2002 (in print).

\endRefs
\vfill
\enddocument

\end